\documentclass[11pt]{amsart}
\usepackage{hyperref,wasysym,cite}
\hypersetup{colorlinks=true, citecolor=BurntOrange, linkcolor=ForestGreen}
\usepackage{amssymb}
\usepackage{graphicx}
\usepackage{amsmath}
\usepackage{caption}
\usepackage{subcaption}
\usepackage[
dvipsnames]{xcolor}
\usepackage[text={6.5in,9in},centering]{geometry}
\usepackage{array}
\newtheorem{Lemma}{Lemma}[section]

\newtheorem{Proposition}[Lemma]{Proposition}

\newtheorem{Theorem}[Lemma]{Theorem}

  {\begin{trivlist}\item[]\textbf{Proof#1 }}%
  {\hspace*{\fill}$\rule{.3\baselineskip}{.35\baselineskip}$\end{trivlist}}

\numberwithin{equation}{section}

\begin{document}

\title{Regime dependent infection propagation fronts in an SIS model}
 
\date{\today}

\begin{abstract}
We show the existence of traveling front solutions in a diffusive classical SIS epidemic model and the SIS model with a saturating incidence in the size of the susceptible population. We investigate the situation where both susceptible and infected populations move around at a comparable rates, but small compared to the spatial scale. In this case, we show that traveling front solutions exist for each fixed positive speed. In the regime where the infected population diffuses slower than the susceptible population, we show the existence of traveling wave solutions for each fixed positive speed and describe their structure and dependence on the wave speed which as it is varied from 0 to infinity. In the regime where the infected population diffuses faster than the susceptible population, we derive a bound for the speeds of the fronts in this regime in which the infection propagates as a front. Moreover, for the classical SIS model we show that there is a case when the spread of the disease is governed by the Burgers-FKPP equation. 
\end{abstract}
\maketitle

\begin{center}    Anna Ghazaryan\textsuperscript{a}\footnote{ Ghazaryan was   funded in part by an AMS-Simons PUI Faculty Research Enhancement Grant.} \footnote{Ghazaryan would like to thank Matt Holzer for a useful discussion  and pointing out the relevant references. The conversation took place during a workshop at Institute for Mathematical and Statistical Innovation (IMSI), which is supported by the National Science Foundation (Grant No. DMS-1929348).},  Vahagn Manukian \textsuperscript{a,b},  Jonathan Waldmann \textsuperscript{a}, Priscilla Yinzime \textsuperscript{a}
 \end{center}
 
 \textsuperscript{a} \address{Department of Mathematics, Miami University,  Oxford, OH 45056, USA,}  
 \par
 \textsuperscript{b} \address{Department of Mathematical and Physical Sciences, Miami University, Hamilton, OH 45011, USA,}  \email{manukive@miamioh.edu.}
 \par 


 
\section{Introduction \label{s:intro}}

In 1927 Kermack and McKendrick published a landmark paper \cite{Kermack} on mathematical modeling in epidemiology. They suggested a compartmental approach to modeling. The whole population is divided into compartments: Susceptible $S$, Infected $I$, and possibly recovered $R$. The Susceptible compartment represents individuals who are susceptible to the disease, the infected compartment represents individuals who are currently infected with the disease and contagious, and the recovered compartment represents individuals who have recovered from the disease and are no longer contagious and have gained immunity. 
The time evolution of these subpopulations may be captured by systems of differential equations.  When only time evolution is considered, these would be ordinary differential equations. In \cite{BrCa, Kuhl, Mart, Schecter} the following general SIS model is described, 
\begin{eqnarray}\label{ode0}
   S_{t}&=&- \beta SI + \gamma I,\\
   I_{t}&=& \beta SI-\gamma I, \nonumber
\end{eqnarray}
where $S$ and $I$, respectively, represent the density of susceptible and infected populations at time $t>0$.  The parameter $\beta>0$ measures the rate at which members of the population are infected and $\gamma>0$ measures the recovery rate of infected members.

A generalization of this model which takes into account the saturating incidence in the size of the susceptible population can be formulated as follows (see for example, \cite{Mart}), 
\begin{eqnarray}\label{ode}
   S_{t}&=&- \frac{\beta SI}{1+\sigma S} + \gamma I,\\
   I_{t}&=& \frac{\beta SI}{1+\sigma S}-\gamma I. \nonumber
\end{eqnarray}
The  saturated incidence rate  in the form  $\frac{1}{1+\sigma S}$ was proposed in \cite{Anderson}.  It captures the behavioral changes in larger susceptible populations that are aimed at preventing infections.  The parameter $\sigma >0$ is called the inhibition constant. The higher values of $\sigma$ correspond to lower incidence rates. 
The system \eqref{ode0}  is the system \eqref{ode} with $\sigma=0$. 

The spatial distribution of the population as well as the random movement of individuals in each subpopulation is often captured by the Laplacian, or, in one spatial dimension, by the second order derivative with respect of the spatial variable. The epidemiological models are then constructed using partial differential equations, more precisely, by systems of reaction-diffusion equations or partly parabolic systems in which some populations diffuse and others do not. For example, after adding  the diffusion terms over the two dimensional space ($(x,y)$-plane), the system \eqref{ode} becomes 
\begin{eqnarray}\label{2d}
    S_{t}&= &d_{1}\Delta S -  \frac{\beta SI}{1+\sigma S} + \gamma I, \\
    I_{t}&= &d_{2}\Delta I + \frac{\beta SI}{1+\sigma S} - \gamma I,\nonumber
\end{eqnarray}
where $\Delta$ is the Laplacian operator with respect to $x$ and $y$ and $S(x,y,t)$ and $I(x,y,t)$ respectively, represent the density of susceptible and infected populations in a location $(x,y)$  at time $t>0$.  The constants $d_{1}$, $d_{2} \geq 0$ represent the rate of diffusion for the susceptible and infected populations, respectively. The subject of this investigation is the existence of traveling waves which are special solutions that preserve their shape while traveling in a preferred direction.
Some initial conditions in \eqref{2d} may evolve into planar waves in this system. Planar waves are waves that propagate along a preferred axis and are constant in the perpendicular direction. Without loss of generality, we may assume that the direction of propagation of the wave is aligned with one of the spatial variables, say, $x$.  We remark that to prove the existence of such waves, it would be sufficient to prove the existence of traveling waves in the system posed on one-dimensional space
\begin{eqnarray} \label{1d}
    S_{t}&= &d_{1}S_{xx} -  \frac{\beta SI}{1+\sigma S} + \gamma I, \\
    I_{t}&= &d_{2}I_{xx} +  \frac{ \beta SI}{1+\sigma S} - \gamma I.\nonumber
\end{eqnarray}
Indeed, if $(S(x,t), I(x,t))$ is a traveling wave solution in \eqref{1d}, then $(S(x,y,t) \equiv S(x,t)$, $I(x,y,t) \equiv I(x,t)$ are the components of the traveling wave solution for the system  \eqref{2d}. 

Therefore, we focus on proving the existence of traveling waves in \eqref{1d}. In particular, we want to prove the existence of traveling fronts, which are traveling waves that asymptotically connect two distinct equilibrium states. 

We investigate \eqref{1d} in the following parameter regimes:
\begin{enumerate}
    \item[Case 1.] $d_1=\alpha d_2 = \mathcal{O}(\epsilon)$, where $\alpha >0$. 
    
\noindent    Case 1 may be interpreted as the situation when the diffusion rates (random movements) of susceptible and infected individuals are comparable, however, they occur at a scale which is much smaller compared to the units in which the spatial variable $x$  is measured.  This may describe a situation where the infected are asymptomatic. 
    \item[Case 2.] $0\le d_2\ll d_1 $.

\noindent     Case 2 assumes that the infection significantly reduces the movement of the infected population, for example, when infected individuals are unable to walk around as much as those who are not infected, or behavior changes in the individual or its social circle occur that prevent the individual from moving around \cite{Lopes}.
    
    \item[Case 3.] $0\le d_1\ll d_2 $. 

\noindent Case 3 assumes that infected individuals diffuse at significantly higher rates. This may be related to infection-induced behavioral changes. For example, the person may become disoriented and venture outside of its regular habitation zone, or a typical ``zombie'' scenario.
\end{enumerate}

Under certain conditions on the parameters in the system, we proved the existence of front solutions for all three cases.   We show that traveling front solutions exist for each fixed positive speed in Cases 1 and 2.   For case 2 we also describe their structure and dependence on the wave speed which as it varies. For Case 3 we were able to derive a condition in the form of the lower bound on  the speed of the propagation of the front. The fronts in this case are phenomenologically different from the fronts in the previous cases and are related to traveling wave solutions of certain nonlinear partial differential equations. 
The obtained results also hold for the case with $\sigma=0$. In the third case, where $\sigma=0$, the spread of the infection is governed by the Burgers-FKPP equation.

The methods used in the paper include applied dynamical systems techniques and Geometric Singular Perturbation theory \cite{Jones, Kuehn}. In the next section we introduce the traveling wave equations and a conserved quantity which the traveling wave equations possess.  The following three sections are devoted to each of the cases described above. In the last section we focus on the situation with  $\sigma=0$ when the system \eqref{1d} reads
\begin{eqnarray} \label{1ds0}
    S_{t}&= &d_{1}S_{xx} -  \beta SI + \gamma I, \\
    I_{t}&= &d_{2}I_{xx} +   \beta SI - \gamma I.\nonumber
\end{eqnarray}
This system is the diffusive version of the classical SIS model. We formulate the implications of our results for this case and demonstrate that the dynamics generated by the traveling wave equations for \eqref{1ds0} is governed by the Burgers-FKPP equation.

\section{Traveling Waves}\label{TW}

Traveling waves are solutions of the underlying partial differential equation that move in a preferred direction while preserving their shape (see, for example, \cite{GLM} and references within). To seek  traveling wave solutions, we introduce a moving coordinate frame $z=x-ct$, where $c$ is the wave speed parameter. Traveling waves are then stationary (time-independent) solutions of the system
\begin{eqnarray} \label{e:sys3}
    S_{t}&=& d_{1}S_{zz} +cS_z-  \frac{\beta SI}{1+\sigma S} + \gamma I,\\
    I_{t}&=& d_{2}I_{zz} + cI_z + \frac{ \beta SI}{1+\sigma S} - \gamma I.\nonumber
\end{eqnarray}
In other words, traveling waves are solutions of the following system of the nonlinearly coupled, ordinary differential equations
\begin{eqnarray}\label{e:sys4}
    0&=& d_{1}S_{zz} +cS_z-  \frac{\beta SI}{1+\sigma S} + \gamma I, \\
   0&=& d_{2}I_{zz} + cI_z+  \frac{\beta SI}{1+\sigma S} - \gamma I. \nonumber
\end{eqnarray}
Traveling fronts are waves that asymptotically connect two different constant states (spatially homogeneous equilibria). There are two lines of constant states in the $(S,I)$-plane for this system. One line is $I=0$ and the other is $S=\gamma/(\beta-\sigma \gamma)$. Not all of these are physical. Since $S$ and $I$ are concentrations of the population, we want  to prove the existence of traveling fronts that satisfy the boundary-like  conditions
\begin{eqnarray}\label{e:bcpde}
    (S,I)(-\infty)=\left(\frac{\gamma}{\beta -\gamma \sigma},1-\frac{\gamma}{\beta -\gamma \sigma}\right), && (S,I)(\infty)=(1,0).
\end{eqnarray}
Note that  parameters $\beta$ and $\gamma$ should satisfy condition 
\begin{equation}
 \beta > \gamma(1+ \sigma)
 \end{equation}
 to guarantee that the constant state at $-\infty$ is physically relevant.

We also observe that the system \eqref{e:sys4} has a conserved quantity which is obtained by adding the equations: 
\begin{equation}\label{e:s1}
d_1U_{zz}+d_2V_{zz}=-cS_z-cI_z.
\end{equation}
Integrating with this expression with the boundary condition at $\infty$, we obtain
\begin{equation}\label{e:fi0}
d_1U_z+d_2I_z+cS+cI=c.
\end{equation}

Traveling fronts are sought as heteroclinic orbits of the first-order system
\begin{eqnarray}
    S_z&=&U,  \label{e:sys5} \\
    d_{1}U_{z}&=&-cU + \frac{\beta SI}{1+\sigma S} - \gamma I, \nonumber\\
    I_z&=&V, \nonumber\\
    d_{2}V_{z}&=&-cV -\frac{\beta SI}{1+\sigma S} + \gamma I, \nonumber
\end{eqnarray}
which is obtained from \eqref{e:sys4} through a coordinate transformation $(S,U, I, V)=(S, S_z, I, I_z)$. We show the existence of solutions that satisfy the boundary conditions \eqref{e:bcpde} by establishing the existence of heteroclinic orbits of the system \eqref{e:sys5} that approaches 
$ (S, U, I, V) =\left(\frac{\gamma}{\beta -\gamma \sigma} ,0,1-\frac{\gamma}{\beta -\gamma \sigma},0\right)$  at  $-\infty$   and  approaches $(S,U,I,V)=(1,0,0,0)$ at $\infty$. 
We call these equilibria
\begin{equation}\label{e:bc1}
    A=\left(\frac{\gamma}{\beta -\gamma \sigma},0,1-\frac{\gamma}{\beta -\gamma \sigma},0\right)  \quad  \mbox{ and } \quad B=(1,0,0,0).
\end{equation}
In this paper consider parameters $\beta$ and $\gamma$ from the region 
 \begin{equation} \label{e:reg}
   R(\gamma, \beta)=\{(\beta, \gamma)\, |  \, \beta > \gamma(1+ \sigma) \}.
 \end{equation}
The conserved quantity \eqref{e:fi0} in the notations of the system \eqref{e:sys5} reads
\begin{equation}\label{e:fi}
d_1U+d_2V+cS+cI=c.
\end{equation}
We use \eqref{e:fi} to write \eqref{e:sys5} as an equivalent lower dimensional system 
\begin{eqnarray}\label{e:vsys2}
d_1 S_z&=&-d_2 V-c\left(I+ S-1\right), \\
I_z&=&V, \nonumber\\
d_2 V_{z}&=&-cV -\frac{\beta SI}{1+\sigma S} + \gamma I.\nonumber
\end{eqnarray}
The equilibria of \eqref{e:vsys2} associated with $A$ and $B$ are then
\begin{equation}\label{e:red_bc1}
    (S,I,V)=\left(\frac{\gamma}{\beta -\gamma \sigma},1-\frac{\gamma}{\beta -\gamma \sigma},0\right), \qquad (S,I,V)=(1,0,0).
\end{equation}

\section{The case of the vanishing diffusion limit.
}

We consider the system \eqref{e:sys3} under the assumption that $d_2 =\epsilon \ll 1$
and $d_1=\alpha d_2$, where $\alpha> 0$ is some constant of order 1, 
\begin{eqnarray} \label{e:sys3_case1}
    S_{t}&=& \alpha \epsilon S_{zz} +cS_z- \frac{\beta SI}{1+\sigma S} + \gamma I,\\
    I_{t}&=& \epsilon I_{zz} + cI_z +\frac{\beta SI}{1+\sigma S}- \gamma I.\nonumber
\end{eqnarray}
The system 
\eqref{e:vsys2} can be written as 
\begin{eqnarray}\label{e:slvansys2}
\epsilon S_z&=&-\frac{\epsilon}{\alpha}  V-\frac{c}{\alpha}\left(I+ S-1\right), \\
I_z&=&V, \nonumber\\
\epsilon V_{z}&=&-cV -\frac{\beta SI}{1+\sigma S} + \gamma I.\nonumber
\end{eqnarray}
Note that \eqref{e:slvansys2} is singularly perturbed system with a small parameter $\epsilon$, which is written in the slow form. It  is associated to the fast system
\begin{eqnarray}\label{e:slvansys3}
S_\xi&=&-\frac{\epsilon}{\alpha}  V-\frac{c}{\alpha}\left(I+ S-1\right), \\
I_\xi&=&\epsilon V, \nonumber\\
V_\xi&=&-cV -\frac{\beta SI}{1+\sigma S} + \gamma I, \nonumber
\end{eqnarray}
through the coordinate transformation $z= \epsilon \xi$.
In the limit as $\epsilon \xrightarrow{} 0 $, the system \eqref{e:slvansys2}
possesses a normally invariant set
\begin{equation}    M_{\epsilon=0}=\left\{(S,I,V)\,|\, I+ S-1=0, -cV -\frac{\beta SI}{1+\sigma S} + \gamma I=0\right\}
\end{equation}
which is  also the set of the equilibrium points of  the system   
\begin{eqnarray}\label{e:slvansys3zero}
S_\xi&=&-\frac{c}{\alpha}\left(I+ S-1\right), \\
I_\xi&=&0, \nonumber\\
V_\xi&=&-cV -\frac{\beta SI}{1+\sigma S} + \gamma I,\nonumber
\end{eqnarray}
which is the system \eqref{e:slvansys3} with $\epsilon=0$.

The linearization of \eqref{e:slvansys3zero} about any point 
$(S,I,V)=\left(1-I, I, 
 -\frac{\beta I(1-I)}{c(1+\sigma (1-I))} + \frac{\gamma}{c} I\right)$
of  $M_{\epsilon=0}$ has two negative  eigenvalues and one zero eigenvalue, and thus $M_{\epsilon=0}$ is normally hyperbolic  and attracting manifold.  By Fenichel's theory \cite{Fenichel2, Fenichel3} there exists a sufficiently small  $\epsilon_0$ such that for any $\epsilon <\epsilon_0$  the system \eqref{e:slvansys2} or equivalently \eqref{e:slvansys3} has a one-dimensional invariant attracting set $M_{\epsilon}$ which is an $\epsilon$-order perturbation of $M_{\epsilon=0}$.  

The reduced flow on $M_{\epsilon=0}$ is given by
\begin{equation}\label{e:red}
 I_z=  \frac{1}{c}   \frac{I\left(\gamma -(\beta- \sigma\gamma)(1-I)\right)}{\beta (1+\sigma (1-I))} .
\end{equation}
The two equilibria of the reduced flow \eqref{e:red} are $I=0$ and $ I= 1- \frac{  \gamma  }{\beta -\sigma\gamma}
$. $I=0$ is a stable equilibrium and $I=1-\frac{\gamma }{\beta -\sigma\gamma}$ is an unstable equilibrium. Since the flow is one-dimensional, there is a connection between $I=1-\frac{  \gamma  }{\beta -\sigma\gamma}$ and $I=0$. 
This is a heteroclinic orbit.  We argue below that this orbit persists for sufficiently small $\epsilon>0$ in \eqref{e:sys3_case1}.

The flow on $M_{\epsilon}$ is given by 
\begin{equation}
    I_z=\frac{1}{c}   \frac{I\left(\gamma -(\beta- \sigma\gamma)(1-I)\right)}{\beta (1+\sigma (1-I))}  +O(\epsilon).
\end{equation}

According to Fenichel theory \cite{Fenichel3}, any invariant set for \eqref{e:slvansys3} that is sufficiently close to $M_{\epsilon=0}$  belongs to $M_{\epsilon}$.
Since for sufficiently small $\epsilon$ the manifold is attracting, the equilibrium $(S,I,V)=\left(\frac{  \gamma  }{\beta -\sigma\gamma} ,1-\frac{  \gamma  }{\beta -\sigma\gamma},0\right)$  must belong to $M_{\epsilon }$ together with its unstable one-dimensional manifold.  Within the one-dimensional set $M_{\epsilon}$, the unstable manifold of $\left(\frac{  \gamma  }{\beta -\sigma\gamma},1-\frac{  \gamma  }{\beta -\sigma\gamma},0\right)$  must intersect the stable manifold of $(S,I,V)=(1,0,0)$ thus forming a heteroclinic orbit in the flow  of \eqref{e:slvansys3}. This heteroclinic orbit is associated with a translationally invariant family of fronts in the original pde. We proved the following proposition.
\begin{Proposition} \label{P:1}
For every fixed  $\alpha$, $\gamma$, $\beta >0$, and $\sigma\geq 0$ 
  that satisfy condition  $\beta > \gamma(1+ \sigma)$, and for every  $c > 0$,  there exists   $\epsilon_0=\epsilon_0(\gamma, \beta,\sigma, c)>0$  such that for each $0<\epsilon <\epsilon_0$  there is a heteroclinic orbit of the system \eqref{e:slvansys2}  that asymptotically connects  $(S,V, I)=\left(\frac{  \gamma  }{\beta -\sigma\gamma},0, 1-\frac{  \gamma  }{\beta -\sigma\gamma}\right)$ at $-\infty$ and $(S,V, I)=(1,0, 0)$ at $+\infty$. 
\end{Proposition}
The main result of this section follows from this proposition. We have the following theorem.
\begin{Theorem} \label{T:1}
 fixed $\alpha$, $\gamma$, $\beta >0$, and $\sigma\geq 0$ 
  that satisfy condition $\beta > \gamma(1+ \sigma)$, and for every  $c > 0$,  there exists   $\epsilon_0=\epsilon_0(\gamma, \beta,c)>0$  such that for each $0<\epsilon <\epsilon_0$  there is an invariant with respect to the translation family of  fronts  of the system \eqref{e:sys3_case1} which move with speed $c$ and have the  rest states  $(S,I)=\left(\frac{  \gamma  }{\beta -\sigma\gamma},1-\frac{  \gamma  }{\beta -\sigma\gamma}\right)$ at $-\infty$ and  $(S,I)=(1,0)$ at $+\infty$.  
\end{Theorem}

\section{Case of decreased mobility of infected population}

In this subsection we assume that the infection slows down the affected population significantly. From a mathematical point of view,  we assume that $d_1=1$ and $ d_2 =\epsilon$, where $0<\epsilon \ll 1$.
The system \eqref{e:sys3} then reads
\begin{eqnarray} \label{e:sys3_case2}
    S_{t}&=&  S_{zz} +cS_z- \frac{\beta SI}{1+\sigma S}  + \gamma I,\\
    I_{t}&=& \epsilon  I_{zz} + cI_z +\frac{\beta SI}{1+\sigma S}  - \gamma I.\nonumber
\end{eqnarray}
The traveling wave equations \eqref{e:vsys2} become
\begin{eqnarray}\label{e:fsys22}
 S_z&=&-\epsilon V-c\left(I+ S-1\right), \\
I_z&=&V, \nonumber\\
\epsilon V_{z}&=&-cV -\frac{\beta SI}{1+\sigma S}  + \gamma I.\nonumber
\end{eqnarray}
Here $V$ is the fast variable,  while $S$, $I$  are the slow variables. We refer to the system \eqref{e:fsys22} as the fast system as opposed to the equivalent
 slow system
\begin{eqnarray}\label{e:slsys22}
 S_{\xi}&=&-\epsilon^2 V-\epsilon c (I+ S-1),\\
I_{\xi}&=&\epsilon V, \nonumber\\
 V_{\xi}&=&-cV -\frac{\beta SI}{1+\sigma S}  + \gamma I,\nonumber
\end{eqnarray}
where  $z= \epsilon \xi$, as in the previous section.
Following the Geometric Singular Theory protocol, we consider the limits of \eqref{e:fsys22} and \eqref{e:slsys22} as $\epsilon \to 0$.  The system \eqref{e:fsys22} with $\epsilon = 0$ gives us the description of a slow manifold 
\begin{eqnarray}
   N_{\epsilon=0}&=&\left\{(S,I,V)\,|\,  V=- \frac{\beta}{c}\left(  \frac{S}{1+\sigma S}  - \frac{\gamma}{\beta}\right) I\right\},
\end{eqnarray}
which is also the 2-dimensional set of equilibria of \eqref{e:fsys22} with $\epsilon = 0$,
\begin{eqnarray}\label{e:slsys22zero}
 S_{\xi}&=&0,\\
I_{\xi}&=&0, \nonumber\\
 V_{\xi}&=&-cV -\frac{\beta SI}{1+\sigma S}  + \gamma I.\nonumber
\end{eqnarray}
The linearization of \eqref{e:slsys22zero} about any point of $N_{\epsilon=0}$  has two zero eigenvalues and one negative eigenvalue $-c$. Therefore, $N_{\epsilon=0}$ is normally hyperbolic and attracting set. 
 By Fenichel's theory \cite{Fenichel2, Fenichel3} there exist sufficiently small  $\epsilon_0$ such that for any $\epsilon <\epsilon_0$  the system \eqref{e:fsys22} or equivalently \eqref{e:slsys22} have an invariant, 2-dimensional, attracting set $N_{\epsilon}$ which is an $\epsilon$-order perturbation of $N_{\epsilon=0}$, 
 \begin{eqnarray}
   N_{\epsilon=0}&=&\left\{(S,I,V)\,| \, V=- \frac{\beta}{c}\left(  \frac{S}{1+\sigma S}  - \frac{\gamma}{\beta}\right) I+O(\epsilon)\right\}.
\end{eqnarray}
 
The reduced flow on $N_{\epsilon=0}$  is given by 
the planar system
\begin{eqnarray}\label{e:redslsys22}
 S_z&=&- c (I+ S-1),\\
I_z&=&-\frac{\beta}{c} I\left( \frac{ S}{1+\sigma S}   - \frac{\gamma}{\beta}\right). \nonumber
\end{eqnarray}
The flow on $N_{\epsilon}$  is an $\epsilon$-order perturbation of \eqref{e:redslsys22},
\begin{eqnarray}\label{e:redslsys22pert}
 S_z&=&- c (I+ S-1),\\
I_z&=&-\frac{\beta}{c} I\left( \frac{S}{1+\sigma S}   - \frac{\gamma}{\beta}\right)+O(\epsilon). \nonumber
\end{eqnarray}

We first focus on the system \eqref{e:redslsys22}. Its  equilibria  are 
\begin{equation}\tilde{B}=(1,0)\,\,\, \mbox{ and }\,\,\, \tilde{A}=\left(\frac{\gamma}{\beta-\gamma\sigma }, 1-\frac{\gamma}{\beta-\gamma\sigma }\right).\end{equation} 
Note that the equilibrium $\tilde A$ is within the first quadrant of the $(S, I)$ plane if $\beta > \gamma(1+ \sigma)$, therefore it is physically relevant when the parameters $(\gamma, \beta)$ are from the region $R(\gamma, \beta)$ defined in \eqref{e:reg}.

The linearization of the vector field generated by \eqref{e:redslsys22} at $\tilde{B}=(1,0)$ has the eigenvalues
\begin{equation} \lambda_1(\tilde{B})=-\frac{\beta -\gamma(1+\sigma)}{c(1+\sigma)},\qquad
    \lambda_2(\tilde{B})=-c.
\end{equation}
Since $c>0$ and $\beta > \gamma(1+ \sigma)$, both eigenvalues $\lambda_1(\tilde{B})$ and $\lambda_2(\tilde{B})$ are negative, therefore the equilibrium $\tilde{B}$ is a stable node.

The linearization of the system \eqref{e:redslsys22} at $\tilde{A}$ has the following eigenvalues
\begin{equation}
    \lambda_1(\tilde{A})=\frac{-c+\sqrt{c^2+4\frac{(\beta-\sigma\gamma)(\beta-(1+\sigma)\gamma)}{\beta}}}{2},\quad
    \lambda_2(\tilde{A})=\frac{-c-\sqrt{c^2+4\frac{(\beta-\sigma\gamma)(\beta-(1+\sigma)\gamma)}{\beta}}}{2}.
\end{equation}
Since $c$ is positive and $\beta > \gamma(1+ \sigma)$, then $\lambda_1(\tilde{A})$ is positive and $\lambda_2(\tilde{A})$ is negative, and, therefore, the equilibrium $\tilde{A}$ is a saddle.

We want to establish the existence of a heteroclinic orbit in the flow generated by \eqref{e:redslsys22} that connects the saddle point $\tilde{A}$ to $\tilde{B}$. We do this by using a trapping region argument and then we apply Fenichel's theory \cite{Fenichel2, Fenichel3} to extend the existence result in the limiting case $\epsilon=0$ to the full system \eqref{e:slsys22} for $0<\epsilon \ll 1$.

To apply the trapping region argument, we analyze the vector field generated by \eqref{e:redslsys22} on the line segments
\begin{eqnarray*}
l_1 &=& \left\{(S,I)\,|\, S=\frac{\gamma}{ \beta - \gamma \sigma}, 0<I<1-\frac{\gamma}{ \beta - \gamma \sigma }\right\},\\
l_2 &=& \left\{ (S,I)\,|\, S=1-I, 0<I<1-\frac{\gamma}{ \beta - \gamma \sigma }\right\},\\
l_3 &=& \left\{(S,I)\,|\, I=0, \frac{\gamma}{ \beta - \gamma \sigma }<S<1\right\},
\end{eqnarray*}
that form the triangle $\triangle \tilde A\tilde B D $. 
Evaluation of the vector field \eqref{e:Fresc_lsys22} at each point $p(S,I)\in  l_1$ yields the vector $$V_1(p)=\left(-I-\frac{\gamma}{ \beta - \gamma \sigma }+1,0\right), \,\,\text{ where } \,\, 0<I<1-\frac{\gamma}{ \beta - \gamma \sigma },$$ therefore the  vector field enters the triangle $\triangle \tilde A\tilde BD $ at each point on $l_1$.

We evaluate the vector field \eqref{e:Fresc_lsys22} at each point $p(S,I)\in l_2$ and obtain the vector \begin{equation} V_2(p)=
\left(0,-\frac{\beta}{c} I\left( \frac{1-I}{1+\sigma (1-I)}  - \frac{\gamma}{\beta}\right) \right), \,\,\text{ where } \,\, 0<I<1-\frac{\gamma}{ \beta - \gamma \sigma }, \end{equation} 
therefore every vector at each point on $l_2$ points inside $\triangle \tilde A\tilde B D $.

\begin{figure}
    \centering
    \includegraphics[width=0.6\linewidth]{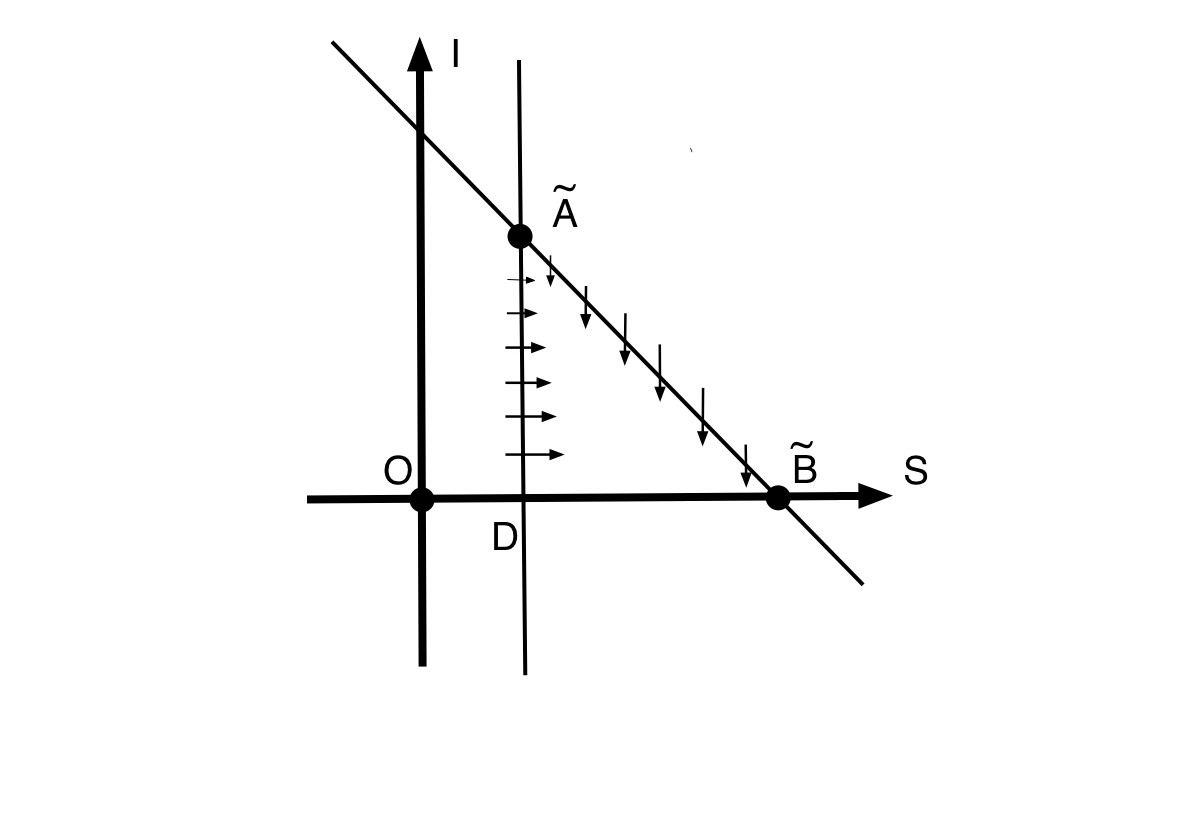}
    \caption{The trapping region for \eqref{e:redslsys22}.}
    \label{fig:1}
\end{figure}

Since the line $I=0$ is invariant for \eqref{e:Fresc_lsys22}, no solution can escape from the triangle $\triangle \tilde A\tilde B D $ through the side $l_3$, so the triangle $\triangle \tilde A\tilde B D $ is a trapping region as shown in Figure~\ref{fig:1}. It is easy to see that the eigenvector that corresponds to the eigenvalue $\lambda_1(\tilde A)$ points inside the triangle. 

\begin{figure}
\begin{subfigure}{.49\textwidth}
  \centering
  \includegraphics[width=1.2\linewidth]{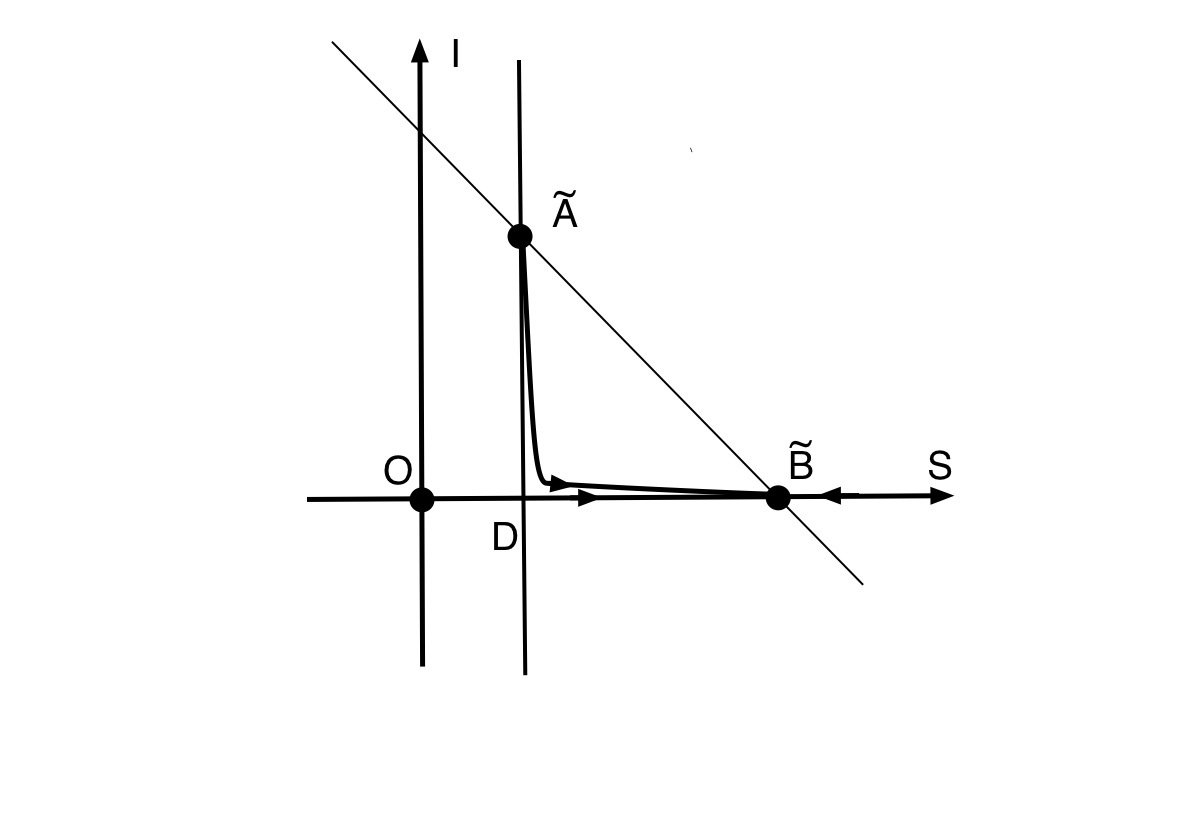}
  \caption{$0<\frac{1}{\delta}\ll 1$ ($c$ is close to zero).}
  \label{fig:sfig1}
\end{subfigure}
\begin{subfigure}{.49\textwidth}
  \centering
  \includegraphics[width=1.2\linewidth]{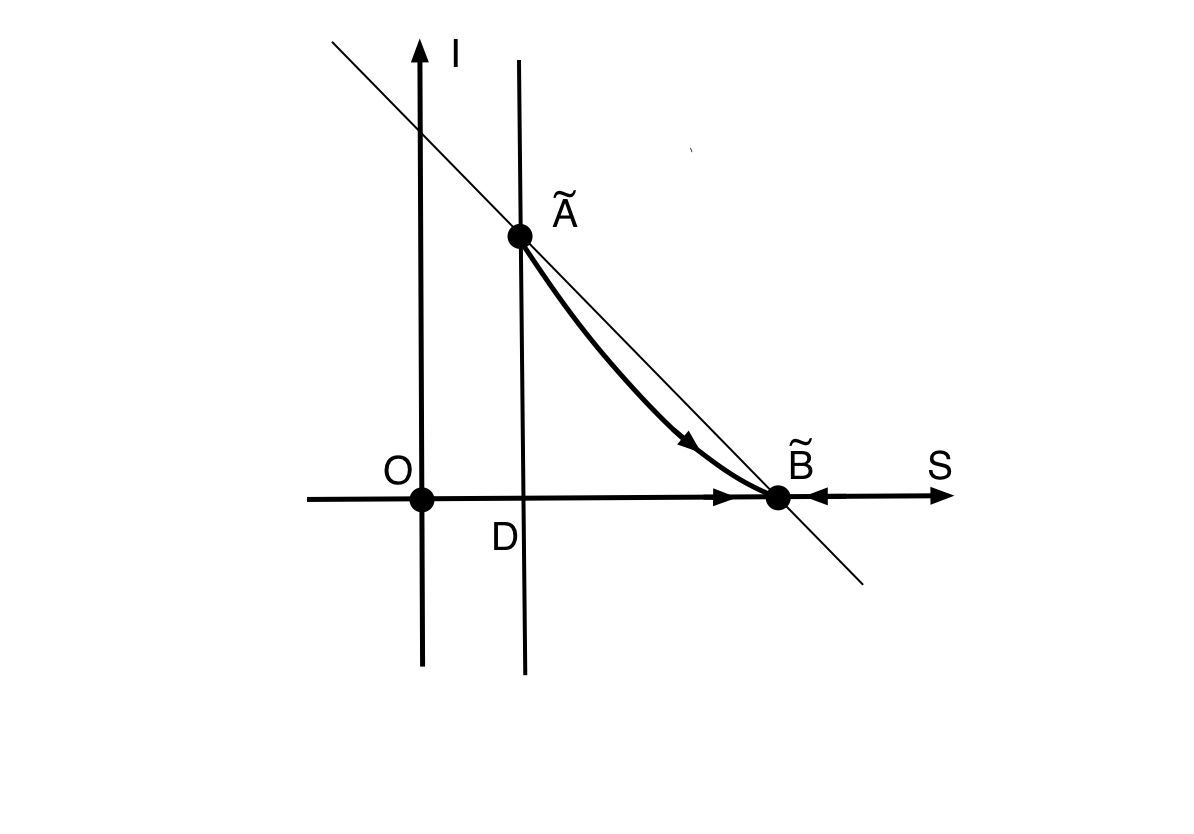}
  \caption{$0<\delta \ll 1$ (large $c$).}
  \label{fig:sfig3}
\end{subfigure}
\caption{Heteroclinic orbits of \eqref{e:resc_lsys22_2} in two singular limits with respect to $\delta$  or, respectively, $c$.}
\label{fig:singrot}
\end{figure}

\begin{figure}
  \centering
  \includegraphics[width=.7\linewidth]{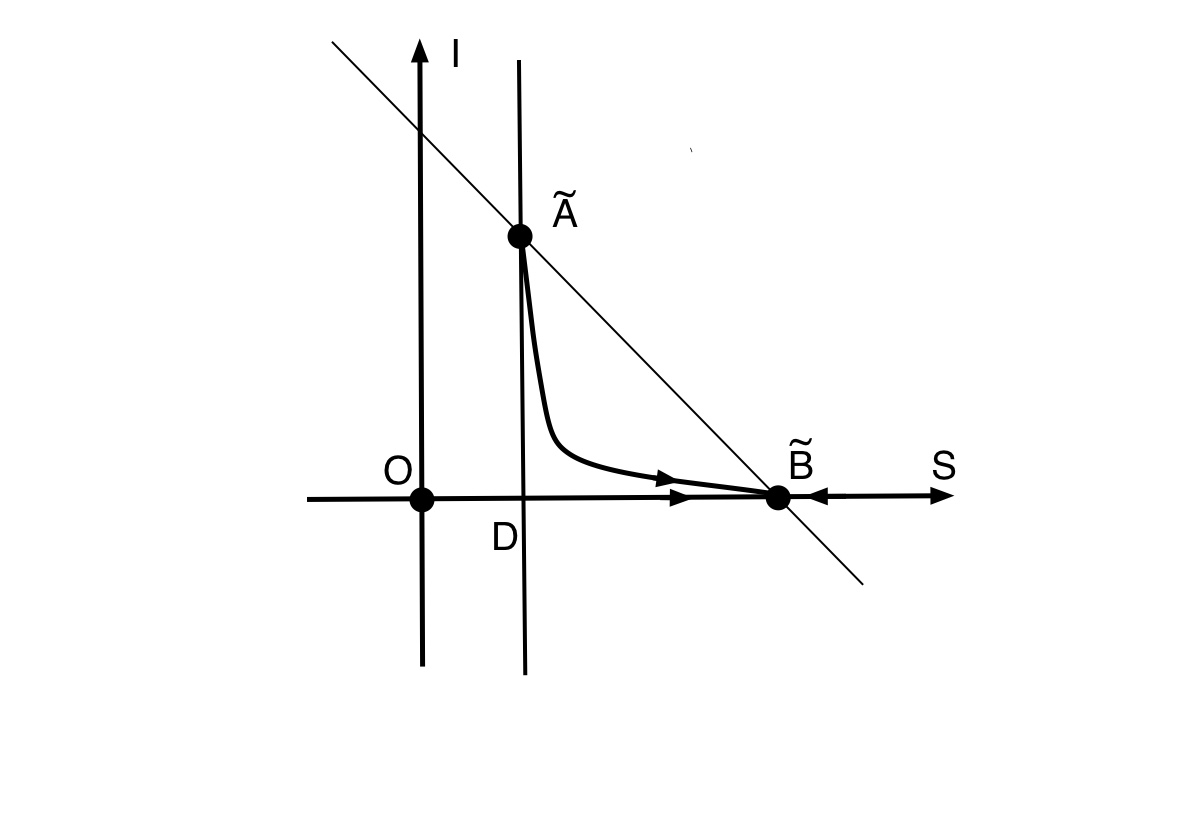}
\caption{Heteroclinic orbit of \eqref{e:resc_lsys22_2} for an intermediate value of $\delta$  (or $c$).}
\label{fig:rot}
\end{figure}

 We use the Poincar\'e-Bendixson theorem to complete the proof of the following lemma: 
\begin{Lemma}\label{L:2} For every fixed $c>0$ and every fixed $\gamma$ and $\beta$ such that $\gamma <\beta$, the system \eqref{e:Fresc_lsys22} (equivalently, \eqref{e:Fresc_lsys22_2}) has heteroclinic orbit that connects $\tilde B = (1, 0)$ and 
$\tilde A =
\left(\frac{\gamma}{ \beta - \gamma \sigma }, 1- \frac{\gamma}{ \beta - \gamma \sigma }\right).$
\end{Lemma}
In addition, we demonstrated how these orbits are related to each other as $c$ changes from $0$ to $\infty$. The details, which are illustrated in Figures~\ref{fig:singrot} and \ref{fig:rot}, are provided later in this section.  

Within the attracting two-dimensional slow manifold $N_{\epsilon=0}$, these orbits are formed by the intersection of the one-dimensional unstable manifold $W^{s}(\tilde{A})$ of the equilibrium $\tilde A$  and the stable two-dimensional manifold $W^{u}(\tilde{B})$ of the equilibrium $\tilde B$. By dimension count, this intersection is transversal and as such persists in the attracting two-dimensional invariant manifold $N_{\epsilon}$ in \eqref{e:slsys22}. We proved the following proposition.
 \begin{Proposition} \label{P:2}
For every fixed  $\sigma\geq 0$,  $\gamma$ and $\beta >0$ that satisfy condition $\beta > \gamma(1+ \sigma)$, and every  $c > 0$,  there exists   $\epsilon_0=\epsilon_0(\gamma, \beta,c)>0$  such that for each $0<\epsilon <\epsilon_0$  there is a heteroclinic orbit of the system \eqref{e:slsys22} that asymptotically connects $(S,V, I)=\left(\frac{\gamma}{ \beta - \gamma \sigma },0, 1-\frac{\gamma}{ \beta - \gamma \sigma }\right)$ as $-\infty$ and $(S,V, I)=(1,0, 0)$ at $+\infty$. 
\end{Proposition}
 
Therefore, we proved the following theorem.

\begin{Theorem} \label{T:2}
For every fixed  $\sigma\geq 0$,  $\gamma$ and $\beta >0$ that satisfy condition $\beta > \gamma(1+ \sigma)$,  and every  $c > 0$,  there exists   $\epsilon_0=\epsilon_0(\gamma, \beta,c)>0$  such that for each $0<\epsilon <\epsilon_0$  there is an invariant with respect to translation family of fronts of the system \eqref{e:sys3_case2} which move with speed $c$ and have the rest states  $(S,I)=\left(\frac{\gamma}{ \beta - \gamma \sigma },1-\frac{\gamma}{ \beta - \gamma \sigma }\right )$ at $-\infty$ and $(S,I)=(1,0)$ at $+\infty$.  
\end{Theorem}

We can say more about the shape of the heteroclinic orbits in the system \eqref{e:redslsys22}. More precisely, we can describe how the shape changes as the traveling velocity parameter $c$ increases from $0$ to $\infty$.

We consider the scaling $\zeta =\frac{z}{c}$ and denote  $\delta =\frac{1}{c^2}$ in \eqref{e:redslsys22},
\begin{eqnarray}\label{e:resc_lsys22_2}
 \delta S_{\zeta}&=&- I- S+1,\\
I_{\zeta}&=&- \beta I\left( \frac{S}{1+\sigma S}  - \frac{\gamma}{\beta}\right). \nonumber
\end{eqnarray}
The associated fast system  in the variable $\eta= \zeta/\delta$ then reads
\begin{eqnarray}\label{e:Fresc_lsys22}
S_{\eta}&=&- I- S+1,\\
I_{\eta}&=&- \delta \beta I\left( \frac{S}{1+\sigma S}  - \frac{\gamma}{\beta}\right). \nonumber
\end{eqnarray}
Note that \eqref{e:Fresc_lsys22} is  a singularly perturbed system when  $\delta$ is  small, equivalently  when the $c$ is large.  In the limit $\delta \rightarrow 0$, the set 
\begin{equation}
\mathcal{L}_0=\{(S,I) \,|\, S=1-I\}
\end{equation}
is an invariant set for \eqref{e:Fresc_lsys22}, which is also a set of equilibria for  the limit of the fast system 
\begin{eqnarray}\label{e:Fresc_lsys22_2}
S_{\eta}&=&- I- S+1,\\
I_{\eta}&=& 0. \nonumber
\end{eqnarray}
The linearization  \eqref{e:Fresc_lsys22_2} about each point of $\mathcal{L}_0$  has a zero eigenvalue and a negative eigenvalue, therefore $\mathcal{L}_0$ is a normally attracting set.
 The reduced flow on $\mathcal{L}_0$ is given by the first order equation
\begin{equation}\label{e:red_lsys22_3}
I_{\zeta}=-  \beta I\left(  \frac{1-I}{1+\sigma (1-I)} -\frac{\gamma}{\beta}\right),
\end{equation}
which has a stable node  at $I=0$ that corresponds to $\tilde{B}$ and a saddle  at $I=1- \frac{\gamma}{ \beta - \gamma \sigma }$ that corresponds to $\tilde{A}$.  
Since $\mathcal{L}_0$ is normally attracting, any portion of $\mathcal{L}_0$, that contains the equilibria $\tilde{A}$ and $\tilde{B}$, for small $\delta\ll 1$ perturbs to a portion of  $\mathcal{L}_\delta$ that contains both equilibria 
$A_\delta$ and $B_\delta$ along with the heteroclinic orbit that asymptotically connects these equilibrium points. 

Now we consider the parameter dependent vector field for \eqref{e:Fresc_lsys22} generated by 
\begin{equation}
  \vec{F}(S,I,\delta)=\begin{pmatrix}
  f_1(S,I)\\
  f_2(S,I)
  \end{pmatrix} = \begin{pmatrix}
 - (I+ S-1)\\
  - \delta I\left(  \frac{S}{1+\sigma S}- \frac{\gamma}{ \beta }\right)
  \end{pmatrix}. 
\end{equation}
For each fixed $0<\frac{\gamma}{ \beta - \gamma \sigma }<1$, we look at the behavior of the vector field  upon changing $\delta$ (or $c$). Since $\frac{\gamma}{ \beta - \gamma \sigma }$ is fixed, the equilibria and nullclines of the system, which play key roles in the organization and the structure of the vector field, are independent of $\delta$. In order to describe the structure of the family of heteroclinc orbits parameterized by $\delta$ or the wave speed, we evaluate the wedge product 
\begin{equation}
\vec{F}(S,I)\wedge \frac{d}{d \delta} \vec{F}(S,I)=(I+ S-1)I\left(  \frac{S}{1+\sigma S}  - \frac{\gamma}{\beta}\right).
\end{equation}
We establish that it is negative at each point inside $\triangle \tilde A\tilde BD $. This wedge product goes into evaluation of the rate of change of the angle $r=\arctan\left(\frac{f_2}{f_1}\right)$ that each vector $\vec{F}(S,I,\delta)$ makes with the $S$-axis as a function of the parameter $\delta$ and 
\begin{equation}
r^\prime(\delta)=\frac{F(S,I,\delta)\wedge \frac{d}{d \delta} F(S,I,\delta)}{|F(S,I,\delta)|^2}.
\end{equation}
Since 
\begin{equation}
F(S,I,\delta)\wedge \frac{d}{d \delta} F(S,I,\delta)=(I+ S-1)I\left(  \frac{S}{1+\sigma S}-\frac{\gamma}{\beta}\right),
\end{equation}
and \begin{equation*}
(I+ S-1)I\left(  \frac{S}{1+\sigma S}-\frac{\gamma}{\beta}\right) <0\end{equation*}
at each interior point of the triangle  $\triangle \tilde A\tilde BD $, so $r^{\prime}(\delta)<0$, and, thus,  the vectors $F(S,I,\delta)$ rotate clockwise monotonically upon changing $\delta$ from $0$ to $\infty$, therefore the saddle separatrix that connects $\tilde A$ to $\tilde B$ rotates from the segment  $\tilde A\tilde B$ when $\delta$ is near 0  to the corner formed by the segments $\tilde AD $ and $D\tilde B$ when $\delta$ is near $\infty$ (or $c$ is close to $0$) as illustrated in Figure~\ref{fig:singrot}A.

\section{The case of small diffusion of the susceptible population}

In this section we study the existence of traveling wave solutions when $\frac{d_1}{d_2}=\epsilon \ll 1$ with $0<\epsilon\ll 1$. Without loss of generality, we take $d_2=1$  and $d_1=\epsilon \ll 1$ and  with $0<\epsilon\ll 1$.  The system \eqref{e:sys3} then reads
\begin{eqnarray} \label{e:sys3_case3}
    S_{t}&=& \epsilon S_{zz} +cS_z- \frac{\beta SI}{1+\sigma S} + \gamma I,\\
    I_{t}&=&   I_{zz} + cI_z +\frac{\beta SI}{1+\sigma S} - \gamma I.\nonumber
\end{eqnarray}
The traveling wave system \eqref{e:vsys2} then can be written 
\begin{eqnarray}\label{e:slsys31}
 \epsilon S_{y}&=&-c (S+ I- 1)-V,\\
I_{y}&=&V, \nonumber\\
 V_{y}&=&-cV -\beta I\left(\frac{S}{1+\sigma S}  - \frac{\gamma}{\beta}\right),\nonumber
\end{eqnarray}
which is singularly perturbed with the small parameter $\epsilon$ within in the slow form. Via the transformation $y=\epsilon\zeta$, \eqref{e:slsys31} can be cast in the equivalent fast form
\begin{eqnarray}\label{e:slsys32}
  S_{y}&=&-c (S+ I- 1)-V,\\
I_{y}&=&\epsilon V, \nonumber\\
 V_{y}&=&\epsilon \left(-cV -\beta I\left(\frac{S}{1+\sigma S}  - \frac{\gamma}{\beta}\right)\right).\nonumber
\end{eqnarray}
Setting $\epsilon=0$ in \eqref{e:slsys31} and \eqref{e:slsys32} we obtain
\begin{eqnarray}\label{e:slsys310}
 0&=&-c (S+ I- 1)-V,\\
I_{y}&=&V, \nonumber\\
 V_{y}&=&-cV -\beta I\left(\frac{S}{1+\sigma S}  - \frac{\gamma}{\beta}\right),\nonumber
\end{eqnarray}
 and 
\begin{eqnarray}\label{e:slsys320}
  S_{y}&=&-c (S+ I- 1)-V,\\
I_{y}&=&0, \nonumber\\
 V_{y}&=&0,\nonumber
\end{eqnarray}
accordingly.
The first equation in  \eqref{e:slsys310} gives us a description  of the slow manifold , which is alos the set of equilibrium points for \eqref{e:slsys320}
\begin{equation}\label{e:manBFKK}
    \mathcal{K}_{\epsilon=0}=\left\{(S,I,V)\,|\,S=1-I-\frac{1}{c} V\right\}.
\end{equation}
Linearizing the system \eqref{e:slsys320} about any point of \eqref{e:manBFKK} we see that the set \eqref{e:manBFKK}
normally hyperbolic and attracting. The reduced flow on the slow manifold \eqref{e:manBFKK} is given by the planar system
\begin{eqnarray}\label{e:pl_sl33}
I_{y}&=& V, \\
 V_{y}&=&
 -cV -\beta I\left(\frac{1-I-\frac{1}{c} V}{(1+\sigma (1-I-\frac{1}{c} V))}  - \frac{\gamma}{\beta}\right).\notag
\end{eqnarray}
The system \eqref{e:pl_sl33} has two equilibria, $\Bar{A}=(1-\frac{\gamma}{\beta-\gamma\sigma},0)$ and  $\Bar{B}=(0,0)$ that respectively correspond to the equilibria $A$ and $B$ of the full system. The linearization of \eqref{e:pl_sl33} at $\Bar{A}$ has the eigenvalues
\begin{eqnarray}\label{ev5}
    \lambda_1(\Bar{A})&=&-c, \\
    \lambda_2(\Bar{A})&=&\frac{\beta-\gamma}{c}+\frac{\gamma  \sigma  \left((1+ \sigma)\gamma -2 \beta  \right) }{c \beta}= \frac{(\beta-(1+\sigma)\gamma)(\beta -\gamma\sigma)}  {c \beta}.\notag
\end{eqnarray}
Since $c>0$ and $\beta> (1+\sigma)\gamma$, then $\lambda_1(\Bar{A})$ is negative and $\lambda_2(\Bar{A})$ is positive. This implies that $\Bar{A}$ is a saddle equilibrium.

The linearization of \eqref{e:pl_sl33} at $\Bar{B}$ has has the two eigenvalues
\begin{eqnarray}
    \lambda_1(\Bar{B})&=&-\frac{c}{2}+\sqrt{\frac{c^2}{4}-\frac{\beta-(1+\sigma) \gamma}{1+\sigma}},\\
    \lambda_2(\Bar{B})&=&-\frac{c}{2}-\sqrt{\frac{c^2}{4}-\frac{\beta-(1+\sigma)  \gamma}{1+\sigma}},\notag
\end{eqnarray}
 which are negative or have negative real parts. Therefore, the equilibrium $\Bar{B}$ is a stable node. Now we will establish the existence of a heteroclinic solution that connects $\Bar{A}$ to $\Bar{B}$ by using a trapping region argument.

We evaluate the vector field $\vec{V}_3$ generated by \eqref{e:pl_sl33} on the following line segments that form a right triangle $\triangle \bar{A}\bar{B}M$: 
\begin{eqnarray*}
s_1 &=& \left\{I,V)\,|\, V=0, \,\,0<I<1-\frac{\gamma}{\beta-\gamma\sigma}\right\},\\
s_2 &=& \left\{ (I,V)\,|\, I=1-\frac{\gamma}{\beta-\gamma\sigma}, \,\,-r\left(1-\frac{\gamma}{\beta-\gamma\sigma}\right)<V<0\right\},\\
s_3 &=& \left\{(I,V)\,|\, V=-rI,\,\, 0<I<1-\frac{\gamma}{\beta-\gamma\sigma}\right\}.
\end{eqnarray*}
Here $r >0$ is the currently undetermined slope of the slanted side of the triangle $\triangle \bar{A}\bar{B}M$.
\begin{figure}
    \centering
    \includegraphics[width=0.6\linewidth]{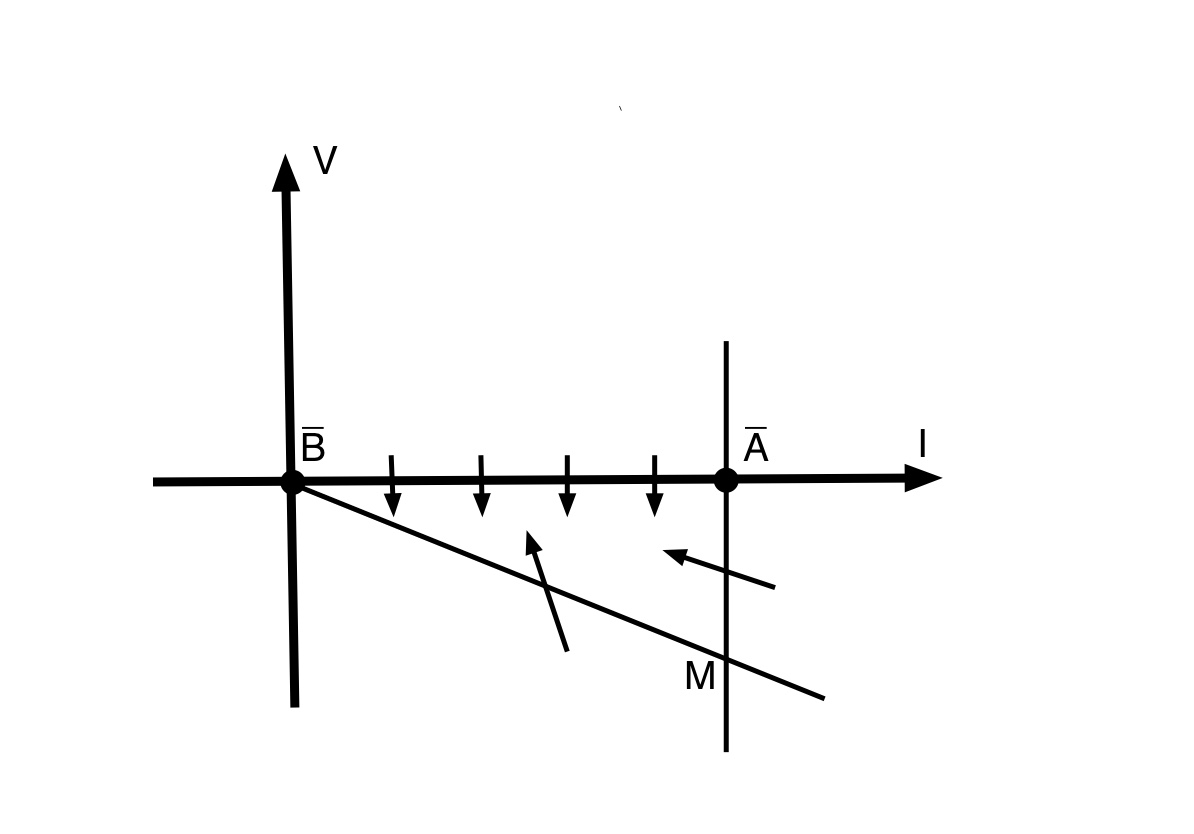}
    \caption{Trapping region for \eqref{e:pl_sl33} }
    \label{fig:2}
\end{figure}
The vector field $\vec{V}_3$ generated by \eqref{e:pl_sl33} on $s_1$ is
\begin{eqnarray}\label{e:s1pl_sl34}
I_{y}&=& 0, \nonumber\\
 V_{y}&=& -\beta I\left(\frac{1-I}{1+\sigma (1-I)}  - \frac{\gamma}{\beta}\right).
\end{eqnarray}

Note that  on $s_1$  $V_{y}$ is negative, since
\begin{equation}\label{e:signs}-\beta I\left(\frac{1-I}{1+\sigma (1-I)}  - \frac{\gamma}{\beta}\right)=-\frac{\beta-\gamma\sigma}{1+\sigma (1-I)} I\left(1-\frac{\gamma}{\beta-\gamma\sigma}-I\right)<0,
    \end{equation}
    where $0<I<1-\frac{\gamma}{\beta-\gamma\alpha}$.
    Therefore, the vector field points down into the triangle $\triangle \bar{A}\bar{B}M $ through the side $\bar{A}\bar{B}$. 

The vector field along the vertical line $s_2$ is given by the system
\begin{eqnarray}\label{e:s2pl_sl35}
I_{y}&=& V, \nonumber\\
 V_{y}&=& 
 -cV -\beta \left(1-\frac{\gamma}{\beta-\gamma\sigma}\right)\left(\frac{ \frac{\gamma}{\beta-\gamma\sigma}-\frac{1}{c} V}{1+\sigma \frac{\gamma}{\beta-\gamma\sigma}-\frac{1}{c} V}  - \frac{\gamma}{\beta}\right). 
\end{eqnarray}
Since $V<0$ on $s_2$, then the vector field points inside $\triangle \bar{A}\bar{B}M$  at each point on the side $\bar{A}M$.

Next, we examine the vector field of \eqref{e:pl_sl33} along the slanted side of the triangle, segment $s_3$.   Parametrically $s_3$ can be described  as 
\begin{equation}
  \begin{pmatrix} I\\ V \end{pmatrix}  =t\begin{pmatrix}
  1-\frac{\gamma}{\beta-\sigma\gamma}\\
  -r\left( 1-\frac{\gamma}{\beta-\sigma\gamma}\right)
  \end{pmatrix} =t\left( 1-\frac{\gamma}{\beta-\sigma\gamma} \right)\begin{pmatrix}
  1\\  -r  \end{pmatrix} ,
\end{equation}
where $t\in(0,1)$. Let us denote $p$ as
$$p= 1-\frac{\gamma}{\beta-\gamma\sigma}. $$  
The vector field evaluated at each point of $s_3$ is 
\begin{equation}
  \vec{V}_3  =\begin{pmatrix}
-rpt\\
tp \left(cr
 -\beta \left(
 \frac{c-ctp+ trp
 }{c+\sigma \left(c-ctp +trp\right)}  - \frac{\gamma}{\beta}\right)\right)
  \end{pmatrix} .
\end{equation}
We calculate the cross product of $\vec{V}_3$ with 
 \begin{eqnarray*}
\vec{d}\times \vec{V}_3 &=&\text{det} \left[\begin{array}{cc}
1 & -r \\ -r p t & tp \left(cr
 -\beta \left(
 \frac{c-ctp+ trp
 }{c+\sigma \left(c-ctp +trp\right)}  - \frac{\gamma}{\beta}\right)\right) 
\end{array}\right]\vec{k}\\
&=&pt\left(-r^2 + cr
 -\beta \left(
 \frac{c-ctp+ trp
 }{c+\sigma \left(c-ctp +trp\right)}  - \frac{\gamma}{\beta}\right)\right) \vec{k},
\end{eqnarray*}
where $ \vec{d} =(1, -r)^T$. 
The vector field points inside the triangle at each point of $s_3$ if the condition 
\begin{equation}\label{e:kpp_cond}
pt\left(-r^2 +cr
 -\beta \left(
 \frac{c-ctp+ trp }{c+\sigma \left(c-ctp +trp\right)}  - \frac{\gamma}{\beta}\right)
 \right)>0
\end{equation}
is satisfied for each $t\in (0,1)$. Since $t\in (0,1)$, the condition \eqref{e:kpp_cond} is equivalent to
\begin{equation}\label{e:kpp_cond3}
r^2 -cr
 +\beta \left(
 \frac{c+tp(r-c)
 }{c+\sigma \left(c+tp(r-c)\right)}  - \frac{\gamma}{\beta}\right)
 <0.
\end{equation}
For the triangle to be a trapping region, we choose $r\leq c$. In fact, if $r> c$, then $c+tp(r-c) > c$. If we think of 
\begin{equation}
    \frac{c+tp(r-c)}{c+\sigma (c+tp(r-c))}
\end{equation}
as $f(c+tp(r-c))$ for $f(x)=\frac{x}{c+\sigma x} $, which is obviously an increasing function and so 
$$f(c+tp(r-c))>f(c)= \frac{1}{1+\sigma},$$
and, therefore,
\begin{equation}\label{e:kpp_cond3_10}
r^2 -cr+\beta \left(
 \frac{c+tp(r-c)
 }{c+\sigma \left(c+tp(r-c)\right)}  - \frac{\gamma}{\beta}\right) > r^2 -cr+ 
 \frac{\beta -(1+\sigma)\gamma 
 }{1+\sigma }.
\end{equation}
However, for $$r>c > \frac{c+\sqrt{c^2-4\frac{\beta-(1+\sigma)\gamma}{1+\sigma}}}{2},$$
due to our assumption on the parameters, we have
\begin{equation} r^2 -cr+ 
 \frac{\beta -(1+\sigma)\gamma 
 }{1+\sigma } > r^2 -cr+ 
 \frac{\beta -(1+\sigma)\gamma 
 }{1+\sigma } >0,\end{equation}
   thus the vector field points outside of the triangle along the side $s_3$. Therefore, we consider only $r \leq c$. In this case $c+tp(r-c) < c$  and   $$f(c+tp(r-c)<f(c)= \frac{1}{1+\sigma},$$ so 
\begin{equation}\label{e:kpp_cond3_1}
r^2 -cr+\beta \left(
 \frac{c+tp(r-c)
 }{c+\sigma \left(c+tp(r-c)\right)}  - \frac{\gamma}{\beta}\right) <  r^2 -cr+ 
 \frac{\beta -(1+\sigma)\gamma 
 }{1+\sigma }.
\end{equation}
The quadratic function on the right has negative values when 
\begin{equation} \label{conditionr}
\frac{c-\sqrt{c^{2}-4  \frac{\beta -(1+\sigma)\gamma 
 }{1+\sigma }}}{2 }<r <\frac{c+\sqrt{c^{2}-4  \frac{\beta -(1+\sigma)\gamma 
 }{1+\sigma }}}{2 }.
 \end{equation} 
This interval is nonempty and contains positive values when
\begin{equation} \label{condonc}
 c\geq 2 \sqrt{ \frac{\beta -(1+\sigma)\gamma 
 }{1+\sigma }}.
\end{equation}
We proved the following lemma.
\begin{Lemma}\label{L:1} For every fixed  $\sigma\geq 0$,  $\gamma$ and $\beta >0$ such that $\beta > \gamma(1+ \sigma)$,  and for any $c$ that satisfies the condition \eqref{condonc},
system \eqref{e:pl_sl33} has a heteroclinic orbit that connects the saddle $\Bar{A}=\left(1-\frac{\gamma}{\beta-\sigma\gamma},0\right)$ and the node $\Bar{B}=(0,0)$.
\end{Lemma}
We point out that the system \eqref{e:pl_sl33}  is equivalent to the second order equation
\begin{equation}\label{e:kpp2}
 I_{yy}+cI_y +\beta I\left(\frac{1-I-\frac{1}{c} I_y}{1+\sigma \left(1-I-\frac{1}{c} I_y\right)}  - \frac{\gamma}{\beta}\right)=0.
\end{equation}
When $\sigma =0$, this equation reduces to the Burgers-FKPP equation and, therefore, more is known about the fronts in this equation and their characteristics \cite{Ma, Mansour, Ryzh}. We discuss this situation in detail in the next section.

Next, we show that the  system \eqref{e:sys3_case3}  supports a front solution that satisfies the condition \eqref{e:bcpde},  
when $0<\epsilon\ll 1$. 
The fronts correspond to heteroclinic orbits in \eqref{e:slsys31} or, equivalently, in \eqref{e:slsys32}.  We argue that heteroclinc orbits of the system \eqref{e:pl_sl33} that we constructed earlier still exist upon switching on the small parameter $\epsilon$. Recall that \eqref{e:pl_sl33} is the reduced flow of \eqref{e:slsys31} or equivalently \eqref{e:slsys32} on the normally attracting manifold \eqref{e:manBFKK}.   Note that the heteroclinic orbits of \eqref{e:pl_sl33} under consideration connect the saddle $\tilde{A}$ to the stable node $\tilde{B}$. These orbits are formed as the intersection of the one-dimensional stable manifold $W^{s}(\tilde{A})$ of the equilibrium $\tilde A$  and the two-dimensional unstable manifold $W^{u}(\tilde{B})$ of the equilibrium $\tilde B$ within a two-dimensional slow manifold. Since  $$\text{dim}\left\{W^{s}(\tilde{A})\right\}+\text{dim}\left\{W^{u}(\tilde{B})\right\}=3> 2,$$  this intersection is transversal. 
 Since the slow manifold $K_{\epsilon=0}$ defined in \eqref{e:manBFKK} normally hyperbolic, by the Fenichel Theory \cite{Fenichel2, Fenichel3}, for any compact subset of \eqref{e:manBFKK}, in the system \eqref{e:slsys31} there exists an invariant two-dimensional, also attracting manifold $K_{\epsilon}$ which is an $\epsilon$-order perturbation of the slow manifold \eqref{e:manBFKK}: $K_{\epsilon} = K_{\epsilon=0}+O(\epsilon)$, more precisely, 
 \begin{equation}\label{e:manBFKKeps}
    \mathcal{K}_{\epsilon}=\left\{(S,I,V)\,|\,S=1-I-\frac{1}{c}V +O(\epsilon)\right\}.
\end{equation} 
This manifold contains the  the equilibria $A$ and $B$. The flow on $\mathcal{K}_{\epsilon}$ is an $\epsilon$-order perturbation of the flow \eqref{e:pl_sl33}
\begin{eqnarray}\label{e:pl_sl35ep}
I_{y}&=& V, \\
 V_{y}&=&
 -cV -\beta I\left(\frac{1-I- \frac{1}{c}V +O(\epsilon)}{1+\sigma (1-I-\frac{1}{c} V)+O(\epsilon)}  - \frac{\gamma}{\beta}\right),\notag
\end{eqnarray}
or
\begin{eqnarray}\label{e:pl_sl35}
I_{y}&=& V, \\
 V_{y}&=&
 -cV -\beta I\left(\frac{1}{\sigma} -\frac{1-I- \frac{1}{c}V }{1+\sigma (1-I-\frac{1}{c} V)}  - \frac{\gamma}{\beta}\right) +O(\epsilon). \notag
\end{eqnarray}
Since the manifold is attracting, the unstable manifold of the equilibrium $B$ must stay on it for all time. Therefore, if a heteroclinic orbit exists it is also confined to this two-dimensional manifold. Within two-dimensional set $K_{\epsilon}$  the  intersection which is transversal  when $\epsilon=0$ (on $K_{\epsilon=0}$) persists for sufficiently small $\epsilon$. Thus, heteroclinic orbits exists in  the system \eqref{e:slsys31} that asymptotically connect $A$ and $B$.
 \begin{Proposition} \label{P:3}
For every fixed  $\sigma\geq 0$,  $\gamma$ and $\beta >0$ such that  $\beta > \gamma(1+ \sigma)$, and for every  $c$  that satisfies the condition \eqref{condonc},
there is   $\epsilon_0=\epsilon_0(\gamma, \beta,c)>0$  such that for each $0<\epsilon <\epsilon_0$  there exists a heteroclinic orbit of the system \eqref{e:slsys31} (equivalently, \eqref{e:slsys31}) that asymptotically connects   $(S,V, I)=\left(\frac{\gamma}{\beta-\sigma\gamma},0, 1-\frac{\gamma}{\beta-\sigma\gamma}\right)$ at $-\infty$ and  $(S,V, I)=(1,0, 0)$ at $+\infty$. 
\end{Proposition}
 
 These orbits in the traveling wave equation represent fronts in the pde system. We proved the following theorem.

\begin{Theorem} \label{T:3}
For every fixed  $\sigma\geq 0$,  $\gamma$ and $\beta >0$  such that  $\beta > \gamma(1+ \sigma)$,     and for every  $c$  that satisfies the condition \eqref{condonc}, there is  $\epsilon_0=\epsilon_0(\gamma, \beta,c)>0$  such that for each $0<\epsilon <\epsilon_0$  there exists an invariant with respect to translation family of fronts of the system \eqref{e:sys3_case3} which move with speed $c$ and have the rest states $(S,I)=\left(\frac{\gamma}{\beta-\sigma\gamma},1-\frac{\gamma}{\beta-\sigma\gamma}\right)$ at $-\infty$ and $(S,I)=(1,0)$ at $+\infty$.  
\end{Theorem}

\section{The case when the inhibition constant is zero.}

Here we discuss the results of the previous sections in the situation when $\sigma=0$.

System \eqref{1ds0} under the assumption that $d_2 =\epsilon \ll 1$ 
and $d_1=\alpha d_2$, with $\alpha>0$, reads
\begin{eqnarray} \label{e:sys3_case1s0}
    S_{t}&=& \alpha \epsilon S_{zz} +cS_z- \beta SI + \gamma I,\\
    I_{t}&=& \epsilon I_{zz} + cI_z +\beta SI - \gamma I.\nonumber
\end{eqnarray}
 The statement of Theorem~\ref{T:1} in this case is as follows. 
\begin{Theorem} \label{T:1s0}
For every fixed $\alpha$, $\gamma$ and $\beta >0$ that satisfy condition $\beta >\gamma$, and for every  $c > 0$,  there exists   $\epsilon_0=\epsilon_0(\gamma, \beta,c)>0$  such that for each $0<\epsilon <\epsilon_0$  there is an invariant with respect to translation family of fronts  of the system \eqref{e:sys3_case1s0} which move with speed $c$ and have the rest states  $(S,I)=\left(\frac{  \gamma  }{\beta },1-\frac{  \gamma  }{\beta }\right)$ at $-\infty$ and $(S,I)=(1,0)$ at $+\infty$.  
\end{Theorem}

In the case when the infection slows down the affected population  ($d_1=1$ and  $ d_2 =\epsilon$, where $0<\epsilon \ll 1$),  system \eqref{1ds0} reads
\begin{eqnarray} \label{e:sys3_case2_s0}
    S_{t}&=&  S_{zz} +cS_z- \beta SI  + \gamma I,\\
    I_{t}&=& \epsilon  I_{zz} + cI_z +\beta SI  - \gamma I,\nonumber
\end{eqnarray}
 and the statement of Theorem~\ref{T:2} is as follows. 
\begin{Theorem} \label{T:2s0}
For every fixed $\gamma$ and $\beta >0$ that satisfy condition $\beta > \gamma$, and every  $c > 0$,  there exists $\epsilon_0=\epsilon_0(\gamma, \beta,c)>0$ such that for each $0<\epsilon <\epsilon_0$  there is an invariant with respect to the translation family of fronts of the system \eqref{e:sys3_case2} which move with speed $c$ and have the rest states  $(S,I)=\left(\frac{\gamma}{ \beta},1-\frac{\gamma}{ \beta}\right )$ at $-\infty$ and  $(S,I)=(1,0)$ at $+\infty$.  
\end{Theorem}

Finally, when $\frac{d_1}{d_2}=\epsilon \ll 1$ with $0<\epsilon\ll 1$   we take  $d_2=1$  and $d_1=\epsilon \ll 1$ in \eqref{1ds0} and obtain
\begin{eqnarray} \label{e:sys3_case3s0}
    S_{t}&=& \epsilon S_{zz} +cS_z- \beta SI + \gamma I,\\
    I_{t}&=&   I_{zz} + cI_z +\beta SI - \gamma I.\nonumber
\end{eqnarray}
The main result in this case is as follows.
\begin{Theorem} \label{T:3s0}
For every fixed  $\gamma$ and $\beta >0$ that satisfy condition $\beta > \gamma$,  and for every  $c$ that satisfies the condition \eqref{condonc}, there is  $\epsilon_0=\epsilon_0(\gamma, \beta,c)>0$  such that for each $0<\epsilon <\epsilon_0$  there exists an invariant with respect to translation family of fronts  of the system \eqref{e:sys3_case3s0} which move with speed $c$ and have the rest states  $(S,I)=\left(\frac{\gamma}{\beta},1-\frac{\gamma}{\beta}\right)$ at $-\infty$ and $(S,I)=(1,0)$ at $+\infty$.  
\end{Theorem}
A noteworthy observation is that in this last case, the equation \eqref{e:kpp2} that governs the dynamics of the system \eqref{e:sys3_case3} is related to the Burgers-FKPP equation. 
Indeed,  when $\sigma=0$ the equation \eqref{e:kpp2} reads
\begin{equation}\label{e:bkpp}
 I_{yy}-\left(\frac{\beta}{c}I-c\right)I_y+\beta I\left(1-\frac{\gamma}{\beta}-I\right)=0.\nonumber
\end{equation}
This is a traveling wave  equation for the following partial differential equation, where we think of $c$ being a ``given,''  fixed parameter
\begin{equation}\label{e:Burgerskpp}
 I_t +\frac{\beta}{c} II_x= I_{xx}+\beta I\left(1-\frac{\gamma}{\beta}-I\right),
\end{equation}
which happens to be the Burgers-FKPP equation \cite{Murray, Xin}
or Fisher-KPP equation with a nonlinear convection \cite{Mansour}. 
In \cite[page 293]{Murray}, the following result is stated: For the equation 
\begin{equation}\label{e:BurgerskppMur}
 u_{\tilde t} +ku u_{\tilde x}= u_{\tilde x \tilde x}+u(1-u)\end{equation}
a traveling front solutions exist for all $\tilde c \geq \tilde c(k)$ where
\begin{equation}\label{Murrayc} \tilde c(k)=\begin{cases} 2 \mbox{ if } k<2,\\
\frac{k}{2}+\frac{2}{k} \mbox{ if } k\geq 2.\end{cases}\end{equation}
The proof of this result can be found in \cite{Ma} where wave speed selection mechanisms are investigated.  

Below we explain that the result in Lemma~\ref{L:1} is aligned with the formula \eqref{Murrayc} \cite[formula (5.73)]{Murray}. 
The equation \eqref{e:Burgerskpp} rescaled as $\tilde I = \frac{\beta}{\beta-\gamma}I$, $\tilde x=\sqrt{\beta -\gamma }\,\,x $, $\tilde t = (\beta -\gamma) t$ is a case  of \eqref{e:BurgerskppMur} with $k= \frac{\sqrt{\beta - \gamma }}{c}$, 
\begin{equation}\label{e:Burgerskpp1}
 \tilde I_{\tilde t} +\frac{\sqrt{\beta-\gamma}}{c} \tilde I\tilde I_{\tilde x}= \tilde I_{\tilde x\tilde x}+ \tilde I(1-\tilde I).
\end{equation}
 Under the assumption \eqref{condonc},  we have $k=\frac{\sqrt{\beta-\gamma}}{c} <\frac{1}{2}$, and, therefore $K <2$. Then according to \eqref{condonc} fronts  exist that move with speed $\tilde c \geq 2$. However, in the original variables $x$ and $t$ in our system, these correspond to the fronts that move with speeds faster than $ 2 \sqrt{\beta -\gamma}$, 
 and so the condition \eqref{Murrayc}  amounts to the condition \eqref{condonc}.
All of these fronts are ``pulled'' fronts \cite{Ryzh}, which means that their long time behavior is determined by the dynamics far ahead of the front interface, rather than by the behavior near the it.

\end{document}